\documentclass[12pt]{article}

\usepackage{amsmath,amssymb,cite}

\def\beq{\begin{equation}}
\def\eqn#1{\beq\label{#1}}
\def\eeq{\end{equation}}

\def\bb {\begin {eqnarray}}
\def\eqnn#1{\bb\label{#1}}
\def\eea {\end {eqnarray}}

\def\bbz{\mathbb{Z}}
\def\bbc{\mathbb{C}}
\def\bac{\bbc} 
\def\bbr{\mathbb{R}}
\def\bbn{\mathbb{N}}

\def\p{{\rm p}}
\def\n{{\rm n}}
\def\u{{\rm \nu}}

\def\nn{\nonumber}

\def\np{\vfill\eject}
\def\nd{\end{document}}

\def\tcc{{\tilde{\cal C}}}

\def\({\left(}
\def\){\right)}
\def\eps{\epsilon}
 
\def\lra{\longrightarrow}

\def\lg{\langle} \def\rg{\rangle} 
\def\vf{\varphi}
\def\ha{{\textstyle{1\over2}}}

\def\gc{{\cal G}^{\bac}}

\def\r{\rho}

\def\bbr{{I\!\!R}}
\def\bbn{I\!\!N}
\def\a{\alpha}
\def\b{\beta}

\def\vr{\vert}

\def\s{\sigma}

\def\L{\Lambda}
\def\bu{$\bullet$}
\def\rank{{\rm rank}}

\def\riga{-\kern-4pt - \kern-4pt -}
\font\fat=cmsy10 scaled\magstep5

\def\Bbullet{\raise-3pt\hbox{\fat\char"0F}}

\def\ca{{\cal A}}  \def\cc{{\cal C}}
\def\cd{{\cal D}}  \def\cf{{\cal F}}
\def\cg{{\cal G}} \def\ch{{\cal H}} 
  
\def\cm{{\cal M}} \def\cn{{\cal N}} 
\def\cp{{\cal P}}  
 \def\ct{{\cal T}}

\def\cih{{\cal C}_\chi}
\def\tcih{{\tilde{\cal C}}_\chi}

\def\ido{intertwining differential operator}
\def\idos{intertwining differential operators}

\def\ca{{\cal A}}
\def\nn{\nonumber}

\def\bu{$\bullet$~~~}

\input epsf.tex
\newcount\figno
\def\fig#1#2#3{
\par\begingroup\parindent=0pt\leftskip=1cm\rightskip=1cm\parindent=0pt
\baselineskip=11pt \global\advance\figno by 1 
\epsfxsize=#3 \centerline{\epsfbox{#2}} \vskip 12pt
#1\par
\endgroup\par}
\def\figlabel#1{\xdef#1{\the\figno}}
\def\encadremath#1{\vbox{\hrule\hbox{\vrule\kern8pt\vbox{\kern8pt
\hbox{$\displaystyle #1$}\kern8pt} \kern8pt\vrule}\hrule}}

\begin{document}

 \textheight=24cm
  \textwidth=16.5cm
  \topmargin=-1.5cm
  \oddsidemargin=-0.25cm

\rightline{preprint SISSA 52/2015/FISI (October 2015)} \vspace{10mm}

\begin{center}

{\LARGE {\bf Classification of Conformal\\[4pt] Representations Induced
from the \\[6pt] Maximal Cuspidal Parabolic}}

\vspace{10mm}

{\bf \large V.K. Dobrev}

\vspace{5mm}

{Scuola Internazionale Superiore di Studi Avanzati}\\
{via Bonomea 265}, {34136 Trieste, Italy}

\vskip 0.5cm {and} \vskip 0.5cm

 Institute for Nuclear
Research and Nuclear Energy\footnote{Permanent address.}\\
{Bulgarian Academy of Sciences}\\ {72
Tsarigradsko Chaussee, 1784 Sofia, Bulgaria}

\end{center}

\vspace{10mm}

\begin{abstract}
In the present paper we continue the project of systematic
construction of invariant differential operators on the example of
representations of the conformal algebra induced   from the maximal
cuspidal parabolic.
\end{abstract}

\vspace{10mm}

\section{Introduction}

Invariant differential operators   play very important role in the
description of physical symmetries. In a recent paper \cite{Dobinv}
we started the systematic explicit construction of invariant
differential operators. We gave an explicit description of the
building blocks, namely, the parabolic subgroups and subalgebras
from which the necessary representations are induced. Thus we have
set the stage for study of different non-compact groups and
induction from different parabolics.

 In the present paper we  focus on the  algebra  ~$so(4,2)$ ~
and representations induced   from the maximal cuspidal parabolic.

 This paper is a   sequel of \cite{Dobinv} and  we refer to
it and \cite{Dobparab}   for motivations and extensive list of literature on the
subject.

\section{Preliminaries}

 Let $G$ be a semisimple non-compact Lie group, and $K$ a
maximal compact subgroup of $G$. Then we have an Iwasawa
decomposition ~$G=KA_0N_0$, where ~$A_0$~ is abelian simply
connected vector subgroup of ~$G$, ~$N_0$~ is a nilpotent simply
connected subgroup of ~$G$~ preserved by the action of ~$A_0$.
Further, let $M_0$ be the centralizer of $A_0$ in $K$. Then the
subgroup ~$P_0 ~=~ M_0 A_0 N_0$~ is a minimal parabolic subgroup of
$G$. A parabolic subgroup ~$P ~=~ M A N$~ is any subgroup of $G$
  which contains a minimal parabolic subgroup.

The importance of the parabolic subgroups comes from the fact that
the representations induced from them generate all (admissible)
irreducible representations of $G$ \cite{Lan,KnZu}.
Actually, induction from the cuspidal parabolic subgroups is enough for this
classification result.

For our purposes here we restrict to ~{\it maximal}~ cuspidal parabolic
subgroups ~$P$, so that $\rank\,A=1$.

Let ~$\nu$~ be a (non-unitary) character of ~$A$, ~$\nu\in\ca^*$,
let ~$\mu$~ fix an discrete series representation ~$D^\mu$~ of ~$M$~ on
a vector space ~$V_\mu\,$, or a limit thereof.

 We consider induced representation ~$\chi =$ Ind$^G_{P}(\mu\otimes\nu \otimes 1)$. These belong
 to the family of
~{\it elementary representations} of $G$ \cite{DMPPT}. (They are called
  {\it generalized principal series representations} (or {\it
limits thereof}) in mathematical literature \cite{Knapp}.)
  Their spaces of functions are:
\eqn{fun} \cc_\chi ~=~ \{ \cf \in C^\infty(G,V_\mu) ~ \vr ~ \cf
(gman) ~=~ e^{-\nu(H)} \cdot D^\mu(m^{-1})\, \cf (g) \} \eeq where
~$a= \exp(H)\in A$, ~$H\in\ca\,$, ~$m\in M$, ~$n\in N$. The
representation action is the $left$ regular action: \eqn{lrr}
(\ct^\chi(g)\cf) (g') ~=~ \cf (g^{-1}g') ~, \quad g,g'\in G\ .\eeq
Note that for our considerations it is enough to use the infinitesimal left action:
\eqn{lac} (X_L\,\cf) (g) ~\doteq~  {d \over dt} \cf
(\exp(-tX)g)\vr_{t=0}\   \eeq

An important ingredient in our considerations are the ~{\it
highest weight representations}~ of ~$\cg$~ associated to the ER ~$\chi\,$. These can be
realized as (factor-modules of) Verma modules ~$V^\L$~ over
~$\cg^\bac$, where ~$\L\in (\ch^\bac)^*$, ~$\ch^\bac$ is a Cartan
subalgebra of ~$\cg^\bac$, the weight ~$\L = \L(\chi)$~ is determined
uniquely from $\chi$ \cite{Dob}. We recall that the
Verma module is explicitly given by ~$V^\L ~=~ U(\gc_-) \otimes v_0\,$,
where we employ the triangular decomposition of ~$\gc$~:
~ $\gc ~=~ \gc_+ \oplus \ch^\bac \oplus \gc_- $,
 ~$\gc_+,\gc_-$~ being the raising,lowering generators
of $\gc$.

 To employ this HWM structure we shall use the {\it right action} of
~$\gc$~   by the standard formula:
\eqn{rac} (X_R\,\cf) (g) ~\doteq~  {d \over dt} \cf
(g\exp(tX))\vr_{t=0}\ ,  \eeq
 where ~$\cf\in \cc_\chi$, ~$g\in G$, first
~$X\in \cg$, then we use complex linear extension to extend this action to
 ~$\gc$.  Note that this action
takes ~$\cf$~ out of ~$\cih$~ for some $X$ but that is exactly why
it is used for the construction of the \idos.

Further we need to
introduce $\bbc$-valued realization $\tcih$ of the space $\cih$ by
the formula:
\eqn{scaz} \vf (g) ~\equiv ~ \lg u_0 , \cf(g) \rg  \eeq
where ~$\lg , \rg$~ is the $M$-invariant scalar product in $V_\mu$, ~$u_0$~ is the
highest weight vector in the discrete series representation  space ~$V_\mu\,$,
or the limit thereof.

On these functions the left/right action of ~$\gc$~ is
defined by:
\eqn{racz} (X_{L/R}\, \vf ) (g) ~\equiv ~\lg u_0 , (X_{L/R}\, \cf )(g) \rg  \eeq

  Generically, Verma modules are irreducible, but for the construction of \idos{} we need
the reducible ones. There is a simple criterion \cite{BGG}:
the {\it Verma module $V^\L$  is reducible} if holds
\eqn{bggr}  {(\L+\r, \b^\vee ) ~=~ m} \ , \quad \b^\vee \equiv 2 \b
/(\b,\b) \ \eeq
where  ~$\b \in\Delta^+$ (the positive roots of ~$(\gc,\ch^\bac)$), $m\in\bbn$,
~$(\cdot\, ,\,\cdot)$~ is the bilinear product in $(\ch^\bac)^*\,$, ~$\rho$~ is half
the sum of the   positive roots.

  Whenever the above  is fulfilled there exists \cite{Dix} in ~$V^\L$~
a submodule which is also a Verma module with shifted weight: ~$V^{\L-m\b}$.
This submodule is generated by a
unique vector $v_s\in V^\L$, called ~{\it singular vector},
such that $v_s\neq \bbc v_0$ and it has the properties
of the highest weight vector of ~$V^{\L-m\b}$~:
\eqnn{lows}
&&  X\, v_s ~=~ (\L - m\b)(X) \cdot v_s ~,
\quad X \in \ch^\bac  \nn\\
&& X\, v_s ~~=~~ 0 ~, \quad X \in \gc_+    \eea

 The above situation will be depicted as follows:
\eqn{subm} V^\L ~\lra ~V^{\L-m\b} \eeq
the arrow points ~{\it to}~ the submodule \cite{Dob}.

 The  singular vector is expressed via \cite{Dob}
 $$  {v^s_{m,\b} = \cp_{m,\b}\, v_0}$$
where ~$\cp_{m,\b}(\cg^\bbc_-)$~ is a polynomial  in the universal
enveloping algebra ~$U(\cg^\bbc_-)\,$.

  Then there exists \cite{Dob}  an ~{\it  \ido{}}~ of order ~$m=m_\b$~:
\eqn{invop}
  \cd_{m,\b} ~:~ \tcc_{\chi(\L)}
~\lra ~ \tcc_{\chi(\L-m\b)} \eeq given explicitly by: \eqn{singvv}
  {\cd_{m,\b} ~=~ \cp_{m,\b}((\cg^\bbc_-)_R) } \eeq where
~$(\cg^\bbc_-)_R$~ denotes the {\it right action} of the elements of ~$\cg^\bbc_-$~ on the functions
~$\vf\,$.

 Thus, in each such situation we have an
 ~{\it invariant differential equation}~ of order ~$m=m_\b$~:
\eqn{invde}  {\cd_{m,\b}\ \vf ~=~ \vf'} \ , \qquad \vf \in \tcc_{\chi(\L)} \ , \quad
\vf' \in \tcc_{\chi(\L-m\b)} \eeq

One main ingredient of our approach is as follows. We group the
(reducible) ERs with the same Casimirs in sets called ~{\it
multiplets} \cite{Dobmul,Dob}. The multiplet corresponding to fixed
values of the Casimirs may be depicted as a connected graph, the
vertices of which correspond to the reducible ERs and the lines
between the vertices correspond to intertwining operators. The
explicit parametrization of the multiplets and of their ERs is
important for understanding of the situation.

In fact, the multiplets contain explicitly all the data necessary to
construct the \idos{}, as we shall demonstrate.

\section{Multiplets of ~{\it $so(4,2)$}~ using maximal cuspidal parabolic}

 Let ~$\cg=so(n,2)$,  $n>2$. It has three nonconjugate parabolic
subalgebras:
\eqnn{parab} \cp ~&=&~ \cm \oplus \ca \oplus \cn \nn\\
\cp_{0} ~&=&~ so(n-2) \oplus \ca_0 \oplus \cn_0 \nn\\
&& \dim\ca_0 =2\ ,\quad \dim\cn_0 = 2(n-1) \nn\\
\cp_{1} ~&=&~ so(n-2) \oplus sl(2,\bbr) \oplus \ca_1 \oplus \cn_1 \\
&& \dim\ca_1 =1\ ,\quad \dim\cn_1 = 2n-3 \nn\\
\cp_{2} ~&=&~ so(n-1,1) \oplus \ca_2 \oplus \cn_2 \nn\\
&& \dim\ca_2 =1\ ,\quad \dim\cn_2 = n \nn\eea
where ~$\cp_0$~ is the minimal parabolic, ~$\cp_1$~ is maximal cuspidal,
~$\cp_2$~ is maximal noncuspidal.  Usually it is the latter that is used, but
in the present paper we study induction from the maximal cuspidal parabolic subalgebra ~$\cp_1\,$.

\subsection{The  ~$so(4,2)$~ main multiplets using maximal cuspidal parabolic}

For  ~$so(4,2)$~  the maximal cuspidal parabolic is:
\eqn{maxc} \cp_1 ~=~ so(2) \oplus sl(2,\bbr) \oplus \ca_1 \oplus \cn_1 \eeq
$ \dim\ca_1 =1$,  ~$\dim\cn_1 = 5$.

The signatures of the ERs are ~
\eqn{sign} \chi_1 ~=~ \{\,n',k,\eps,\nu'\,\}\eeq
 where ~$n'\in\bbz$~ is a character of ~$SO(2)$, ~$\nu'\in\bbc$~ is a
character of ~$A_1\,$, ~$k,\eps$~ fix a discrete series
representation of ~$SL(2,\bbr)$, ~$k\in\bbn$, ~$\eps =\pm1$, or a
limit thereof when $k=0$.

The relation with the ~$sl(4)$~ Dynkin labels is as follows
\cite{DoJMP}: \eqn{intp} m_1 ~=~ \ha (k-\nu'+n') \ ,\quad m_2 ~=~ -k
\ , \quad m_3 ~=~ \ha (k-\nu' -n') \eeq

The main multiplet of reducible ERs contains 12 members which we parametrize
as part of the main ~$sl(4)$~ multiplet with 24 members. Thus, the 12-plet has
the following signatures:
\eqnn{weylsl4p1}
&&\L_2 ~=~ (m_{12},-m_2,m_{23}) ~=~ \{n'~=~m_1-m_3 , k~=~m_2,\nu'~=~-m_{13}\}  \nn\\
&&\L_{12} ~=~ (m_2,-m_{12},m_{13}) ~=~ \{-m_1-m_3 ,m_{12},-m_{23} \}
\nn\\
&&\L_{32} ~=~ (m_{13},-m_{23},m_2) ~=~ \{ m_1+m_3, m_{23},-m_{12} \}
\nn\\
&&\L_{121} ~=~ (-m_2,-m_{1},m_{13}) ~=~ \{-m_{13}-m_2 , m_{1} ,-m_3\}
\nn\\
&&\L_{132} ~=~ (m_{23},-m_{13},m_{12}) ~=~ \{m_3-m_1 , m_{13}, -m_2 \}
\nn\\
&&\L_{232} ~=~ (m_{13},-m_3,-m_{2}) ~=~ \{m_{13}+m_2, m_3,-m_1\}
\nn\\
&&\L_{1232} ~=~ (m_{23},-m_{3},-m_{12}) ~=~ \{m_{13}+m_2 , m_{3}, m_1 \}
\nn\\
&&\L_{1321} ~=~ (-m_{23},-m_{1},m_{12})~=~ \{ -m_{13}-m_2 ,m_{1},m_3\}
 \nn\\
&&\L_{2132} ~=~ (m_{3},-m_{13},m_{1}) ~=~ \{ m_3-m_1 , m_{13},m_2\}
\nn\\
&&\L_{12132} ~=~ (m_3,-m_{23},-m_{1}) ~=~ \{ m_1+m_3,m_{23},m_{12} \}
\nn\\
&&\L_{21321} ~=~ (-m_{3},-m_{12},m_{1}) ~=~ \{-m_1-m_3,m_{12},m_{23} \}
\nn\\
&&\L_{121321} ~=~ (-m_3,-m_{2},-m_{1}) ~=~ \{ m_1-m_3,m_{2} ,m_{13} \}
\eea
where ~
$$\L_{i_1i_2\ldots i_t} ~=~ \s_{i_t}\cdots \s_{i_2}\s_{i_1}\L_0 $$
$\s_j$, ~$j=1,2,3$, are the three simple $sl(4)$ reflections,
~$\L_0$~ is the ER, or rather the corresponding Verma module with dominant
highest weight with signature ~$(m_1,m_2,m_3)$, ~$m_k\in\bbn$, which module
is fixing the  $sl(4)$ 24-plet. We have given the signature in both the $sl(4)$ signature notation
$(\cdot,\cdot,\cdot)$
and in the ~$\cp_1$-induced notation \eqref{sign}.

We would like to follow connections with ERs induced from the maximal noncuspidal
parabolic $\cp_2\,$. Thus, below we shall replace the ~$m_k$~ notation with equivalent one:
$$(p,\nu,n) = (m_1,m_2,m_3)$$

Thus, we give the same 12-plet in ~$p,\nu,n$~ parametrization and adding
the Harish-Chandra (HC) parameters \cite{HC}
for the three nonsimple roots in the order ~$\a_{12},\a_{23},\a_{13}$.
  Thus, in the 4-th,5-th,6-th, place the parameters
are ~:~$m_{12}=m_1+m_2,m_{23}=m_2+m_3,m_{13}=m_1+m_2+m_3$, and of course they are redundant but
some representation theoretic statements are formulated easier in their terms. In particular, the
$K$-noncompact HC  parameters are: ~$m_2,m_{12},m_{23},m_{13}$.

Thus, the 12-plet is given now as:
\eqnn{weylsl4pp11}
&&\L^-_0 = \L_2 = (p+\nu ,-\nu,n+\nu \,;\, p, n, p+\nu+n)  ~=~ \chi'^-_{p\nu n}
 \nn\\
&&\L^-_a =\L_{12} = (\nu,-p-\nu ,p+\nu+n \,;\, -p,n, n+\nu )  ~=~ \chi''^-_{p\nu n}  \nn\\
&&\L^-_b =\L_{32} = (p+\nu+n  ,-n-\nu ,\nu\,;\, p, -n, p+\nu)   ~=~ \chi''^+_{p\nu n} \nn\\
&&\L^-_c =\L_{121} = (-\nu,-p,p+\nu+n \,;\, -p-\nu, n+\nu, n )  \nn\\
&&\L^-_d =\L_{132} = (n+\nu ,-(p+\nu+n)  ,p+\nu \,;\, -p, -n, \nu )  ~=~ \chi'^+_{p\nu n}  \nn\\
&&\L^-_e =\L_{232} = (p+\nu+n  ,-n,-\nu\,;\,  p+\nu, -n-\nu, p)   \nn\\
&&\L^+_e =\L_{1232} = (n+\nu ,-n ,-p-\nu \,;\, \nu, -p-\nu-n, -p)  \nn\\
&&\L^+_d =\L_{2132} = (n ,-(p+\nu+n)  ,p\,;\, -p-\nu, -\nu-n,-\nu)   ~=~ \chi^+_{p\nu n}  \nn\\
&&\L^+_c =\L_{1321} = (-n-\nu ,-p,p+\nu\,;\, -p-\nu-n, \nu,-n)  \nn\\
&&\L^+_b =\L_{12132} = (n,-n-\nu ,-p\,;\, -\nu,-p-\nu-n,-p-\nu )  \nn\\
&&\L^+_a = \L_{21321} = (-n ,-p-\nu ,p\,;\, -p-\nu-n,-\nu, -\nu-n,)  \nn\\
&&\L^+_0 =\L_{121321} = (-n,-\nu,-p\,;\, -\nu-n,-p-\nu, -p-\nu-n)
\eea
where we have also indicated the five cases ~$\chi^{(','')\pm}_{p\nu n}$~ coinciding by signatures
with ERs induced from the maximal noncuspidal parabolic. The notations ~$\L^\pm_0\,$,
~$\L^\pm_a\,$, etc, are used in Fig.1. where we present this 12-plet. The arrows
denote both the \idos\ between ERs and embeddings between  Verma modules.
Notation at the arrows denote  the representation parameter and the root, e.g.,
~$\nu_{1}$~ denotes singular vector (embedding) of weight ~$\nu\,\a_{1}$,
~$\p_{12}$~ denotes singular vector (embedding) of weight ~$p\,\a_{12}$.
Note that only the $M$-non-compact roots are involved \cite{Dob}. Thus, the
$M$-compact root ~$\a_2$~ is not relevant for the \idos.

Some remarks on the ERs content: As in the general situation the ERs
~$\L^+_d~=~ \chi^+_{p\nu n}$~ in \eqref{weylsl4pp11} contain
 holomorphic discrete series representations when the discrete parameter ~$\eps=1$,
 and the antiholomorphic discrete series representations when the discrete parameter ~$\eps=-1$.
 (The criterion for holomorphicity is that the $K$-compact HC parameters are
 positive, while the $K$-noncompact HC parameters are negative \cite{Knapp}.)\\
In the ERs ~$\L^-_a~=~ \chi''^-_{p11}$~ and ~$\L^-_b~=~ \chi''^+_{11n}$~  are contained the massless
representations with conformal weight ~$d=1+j$~ with spin ~$j\geq 1$, where ~$j=(p+1)/2$ and $j=(n+1)/2$, resp.
(the three with  lower spin are considered below).

Our diagrams account also for the Knapp-Stein (KS) \cite{KnSt} integral operator relevant
for ~$P_1$~ induction. It acts on the signatures
as the highest $sl(4)$ root $\a_{13}$, \cite{DoJMP,DoMo}. Thus, on the $P_1$-signature
it acts by changing the
sign of the last entry in the ~$\{ \}$~-notation ($\nu'$).
On the figure the KS operators intertwine the ERs symmetric w.r.t. the dashed line. Of course,
the KS from ~$\L^-_c,\L^-_d,\L^-_e$~ to ~$\L^+_c,\L^+_d,\L^+_e$, resp., degenerated to
differential operators of degrees ~$n,\nu,p$, resp., as shown on Fig.1. The KS operators
from ~$\L^+_c,\L^+_d,\L^+_e$~ to ~$\L^-_c,\L^-_d,\L^-_e$, resp., remain integral operators.

The same remarks about the KS integral operators will be true verbatim
for all further Figures below and will not be repeated.

\subsection{Reduced multiplets}

We have several types of reduced multiplets.

\subsubsection{Symmetrically reduced multiplets}

\bu The first case is a septuplet depending on two parameters, and the signatures
may be obtained by setting formally ~$\nu=0$~ in  \eqref{weylsl4p1}~:
\eqnn{weylsl4p1rnu}
&&\L_2 ~=~ (p,0,n) ~=~ \{n'~=~p-n , k~=~0,\nu'~=~-p-n\} ~=~ \chi'^-_{p0n}~=~ _2\chi^-_{pn}\nn\\
&&\L_{12} ~=~ (0,-p,p+n) ~=~ \{-p-n ,p,-n \} ~=~ \chi''^-_{p0n}
\nn\\
&&\L_{32} ~=~ (p+n,-n,0) ~=~ \{ p+n, n,-p \}  ~=~ \chi''^+_{p0n}
\nn\\
&&\L_{132} ~=~ (n,-p-n,p) ~=~ \{n-p , p+n, 0 \} ~=~ \chi'^+_{p0n}  ~=~ \chi^+_{p0n} ~=~ _2\chi^+_{pn}
\nn\\
 &&\L_{1232} ~=~ (n,-n,-p) ~=~ \{p+n , n, p \}
\nn\\
&&\L_{1321} ~=~ (-n,-p,p)~=~ \{ -p-n ,p,n\}
 \nn\\
   &&\L_{121321} ~=~ (-n,0,-p) ~=~ \{ p-n,0 ,p+n \}
\eea
or adding the HC parameters:
\eqnn{weylsl4p1rnuhc}
&&\L^-_{0} ~=~ \L_2 ~=~ (p,0,n \,;\, p,n,p+n) ~=~   _2\chi^-_{pn} \nn\\
&&\L^-_a ~=~ \L_{12} ~=~ (0,-p,p+n \,;\, -p,n,n)
\nn\\
&&\L^-_b ~=~ \L_{32} ~=~ (p+n,-n,0 \,;\, p,-n,p)
\nn\\
&&\L_c ~=~ \L_{132} ~=~ (n,-p-n,p \,;\, -p, -n, 0) ~=~ _2\chi^+_{pn}
\nn\\
 &&\L^+_b ~=~ \L_{1232} ~=~ (n,-n,-p \,;\, 0, -n-p, -p)
\nn\\
&&\L^+_a ~=~ \L_{3212} ~=~ (-n,-p,p \,;\, -n-p, 0, -n)
 \nn\\
   &&\L^+_0 ~=~ \L_{213213} ~=~ (-n,0,-p \,;\, -n,-p,-n-p)
\eea

The first and fourth entries ~$_2\chi^\pm_{pn}$~ form a doublet w.r.t.
induction from the maximal noncuspidal parabolic.

The first entry  ~$_2\chi^-_{pn}$~ is induced from a limit of ~$sl(2,\bbr)$~ discrete series.
For ~$p=n=1$~ it contains the scalar massless  representation  with conformal weight ~$d=1$.

The fourth entry ~$_2\chi^+_{pn}$~ contains  limits  of
holomorphic/antiholomorphic discrete series (for ~$\eps=\pm1$).

This septuplet is shown on Fig.2.

\subsubsection{Asymmetrically reduced multiplets}

\bu The second  case is also septuplet depending on two parameters,   the signatures
may be obtained by setting formally ~$p=0$~ in  \eqref{weylsl4p1}~:
\eqnn{weylsl4pp11p}
&&\L^-_{0\n} ~=~ \L_2 = (\nu ,-\nu,n+\nu \,;\, 0, n, \nu+n)  ~=~ \chi'^-_{0\nu n}
~=~ \chi''^-_{0\nu n}~=~ _1\chi^-_{\nu n}\nn\\
&&\L^-_b ~=~ \L_{32} = (\nu+n  ,-n-\nu ,\nu\,;\, 0, -n, \nu)   ~=~ \chi''^+_{0\nu n}
~=~ \chi'^+_{0\nu n}~=~ _1\chi^+_{\nu n} \nn\\
&&\L^-_d ~=~ \L_{212} = (-\nu,0,\nu+n \,;\, -\nu, n+\nu, n )  \nn\\
&&\L_e ~=~ \L_{232} = (\nu+n  ,-n,-\nu\,;\,  \nu, -n-\nu, 0)   \nn\\
&&\L^+_d ~=~ \L_{3212} = (-n-\nu ,0,\nu\,;\, -\nu-n, \nu,-n)  \nn\\
&&\L^+_b ~=~ \L_{2132} = (n ,-(\nu+n)  ,0\,;\, -\nu, -\nu-n,-\nu)   ~=~ \chi^+_{0\nu n}
   \nn\\
&&\L^+_{0\n} ~=~ \L_{23212} = (-n ,-\nu ,0\,;\, -\nu-n,-\nu, -\nu-n)
\eea
The first two entries ~$_1\chi^\pm_{\nu n}$~ form a doublet w.r.t.
induction from the maximal noncuspidal parabolic.
The ER ~$_1\chi^\pm_{11}$~  contains one (of the two)  spin 1/2
  massless  representations  with conformal weight ~$d=3/2$.

The third and fifth entries are induced from limits of discrete series.

This septuplet is shown on Fig.3n.

The conjugate septuplet may be obtained by setting $n=0$:
\eqnn{weylsl4pp11n}
&&\L^-_{0\p} = \L_2 = (p+\nu ,-\nu,\nu \,;\, p, 0, p+\nu)  ~=~ \chi'^-_{p\nu 0} ~=~
\chi''^+_{p\nu 0} ~=~ _3\chi^-_{p\nu} \nn\\
&&\L^-_a =\L_{12} = (\nu,-p-\nu ,p+\nu \,;\, -p,0, \nu )  ~=~ \chi''^-_{p\nu 0}
~=~ \chi'^+_{p\nu 0} ~=~ _3\chi^+_{p\nu} \nn\\
&&\L^-_e =\L_{232} = (p+\nu  ,0,-\nu\,;\,  p+\nu, -\nu, p)   \nn\\
&&\L_d =\L_{212} = (-\nu,-p,p+\nu \,;\, -p-\nu, \nu, 0 ) \nn\\
&&\L^+_e =\L_{1232} = (\nu ,0 ,-p-\nu \,;\, \nu, -p-\nu, -p)  \nn\\
  &&\L^+_a = \L_{23212} = (0 ,-p-\nu ,p\,;\, -p-\nu,-\nu, -\nu,) ~=~   \chi^+_{p\nu 0} \nn\\
&&\L^+_{0\p} =\L_{213213} = (0,-\nu,-p\,;\, -\nu,-p-\nu, -p-\nu)
\eea
The interpretation of the members is exactly as for the last case, e.g.,
the first two entries ~$_3\chi^\pm_{p\nu }$~ form a doublet w.r.t.
induction from the maximal noncuspidal parabolic.
The ER ~$_3\chi^\pm_{11}$~  contains the other  spin 1/2
  massless  representation with conformal weight ~$d=3/2$.

This septuplet is shown on Fig.3p.

\subsubsection{Symmetrically doubly  reduced multiplets}

\bu The next case is a quartet depending on  one parameter with signatures
may be obtained by setting formally ~$p=n=0$~ in  \eqref{weylsl4p1}~:
\eqnn{weylsl4pp11pn}
&&\L^-_0 ~=~ \L_2 = (\nu ,-\nu, \nu \,;\, 0, 0, \nu)  ~=~ \chi'^-_{0\nu 0}
~=~ \chi''^-_{0\nu 0}    ~=~ \chi''^+_{0\nu 0}
~=~ \chi'^+_{0\nu 0} ~=~ \chi^s_{\nu} \nn\\
&&\L_d ~=~ \L_{212} = (-\nu,0,\nu \,;\, -\nu, \nu, 0 )  \nn\\
&&\L_e ~=~ \L_{232} = (\nu  ,0,-\nu\,;\,  \nu, -\nu, 0)   \nn\\
&&\L^+_0 ~=~ \L_{2132} = (0 ,-\nu  ,0\,;\, -\nu, -\nu,-\nu)   ~=~ \chi^+_{0\nu 0}
    \eea
  The first entry is a singlet w.r.t. induction from the maximal noncuspidal parabolic.
  For ~$\nu=1$~ it is a Lorentz scalar positive energy representation with conformal weight ~$d=2$, i.e.,
   above the unitarity threshold ~$d=1$, but below the limit of
   holomorphic discrete series  ~$d=3$.
  The second and third entries are induced from limits of discrete series.

This  quartet is shown on Fig.4.

\subsubsection{Asymmetrically doubly reduced multiplets}

\bu The next two conjugate cases are triplets depending on  one parameter. The first
 be obtained by setting formally ~$p=\nu=0$~ in  \eqref{weylsl4p1}~:
\eqnn{weylsl4pp1nu}
&&\L^-_{0\p\u} = \L_2 = (0 ,0,n \,;\, 0, n, n)  ~=~ \chi'^-_{00n}
~=~ \chi''^-_{00n}
\nn\\
&&\L_b = \L_{32} = (n  ,-n ,0\,;\, 0, -n, 0)   ~=~ \chi''^+_{00n}
~=~ \chi'^+_{00n} ~=~ \chi^+_{00n}   \nn\\
  &&\L^-_{0\p\u} ~=~ \L_{3212} = (-n ,0,0\,;\, -n, 0,-n)
\eea
The first and third entries are induced from limits of discrete series.

\bu The conjugate case:
\eqnn{weylsl4p1rnuhcnnu}
&&\L^-_{0\n\u} ~=~ \L_2 ~=~ (p,0,0 \,;\, p,0,p) ~=~  \chi'^-_{p00}  ~=~ \chi''^+_{p00} \nn\\
&&\L_a ~=~ \L_{12} ~=~ (0,-p,p \,;\, -p,0,0)   ~=~ \chi''^-_{p00} ~=~ \chi'^+_{p00}  ~=~ \chi^+_{p00}
\nn\\
     &&\L^+_{0\n\u} ~=~ \L_{121321} ~=~ (0,0,-p \,;\, 0,-p,-p)
\eea

The two conjugate triplets are shown on Fig. 5p,5n.


\vspace{10mm}


\section*{Acknowledgments}
The author would like to thank for hospitality the International
School for Advanced Studies (SISSA), Trieste, where part of the work was
done.\\    The author has received partial support
    from European COST actions MP-1210 and MP-1405,
  and from the Bulgarian National Science Fund Grant DFNI-T02/6.


\vspace{10mm}

\fig{}{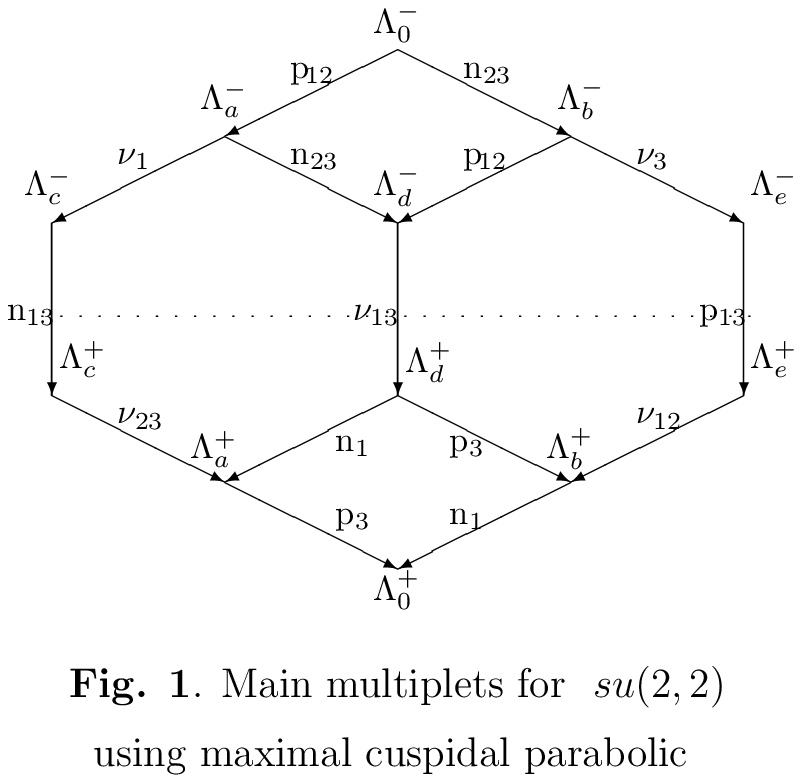}{14cm}

\np

\fig{}{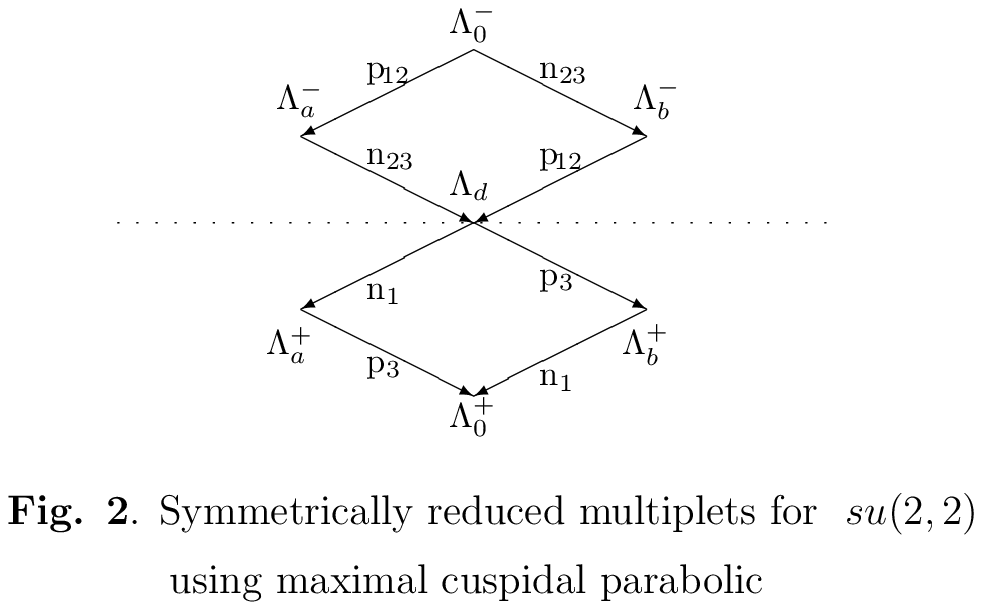}{12cm}

\vspace{10mm}

\fig{}{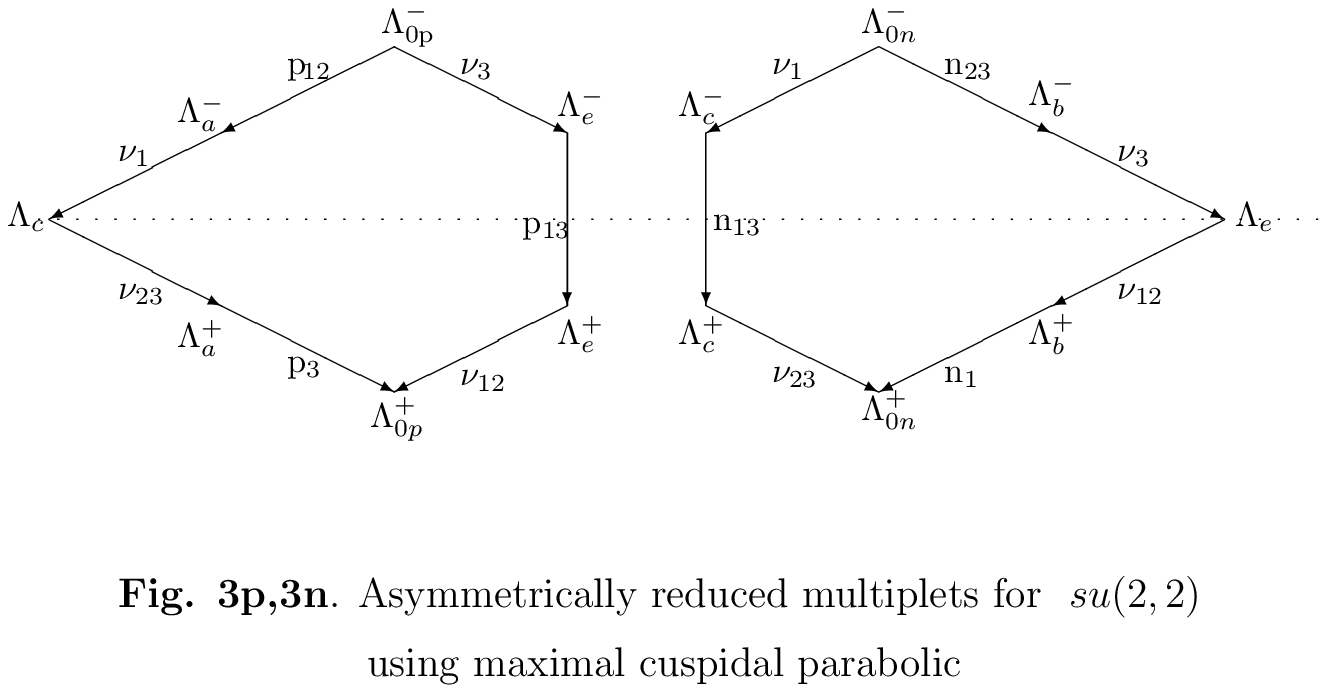}{16cm}

\np

\fig{}{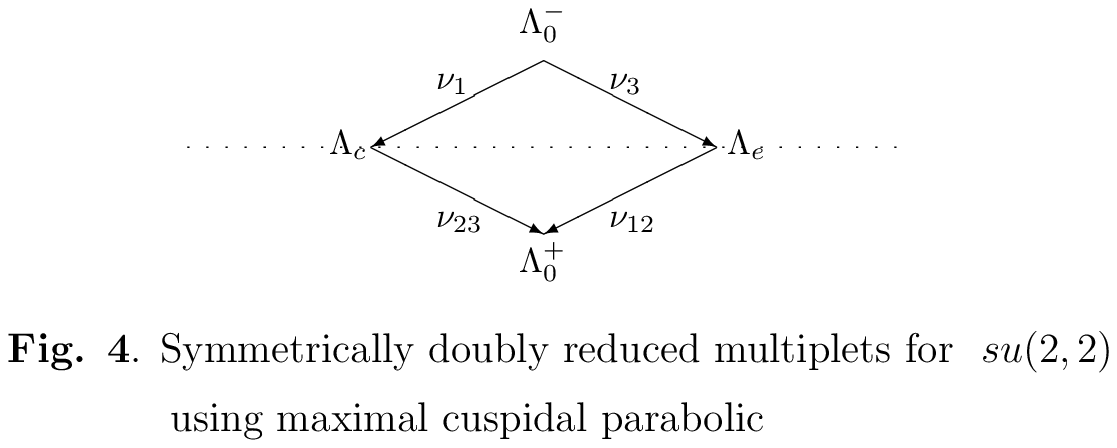}{16cm}

\vspace{10mm}

\fig{}{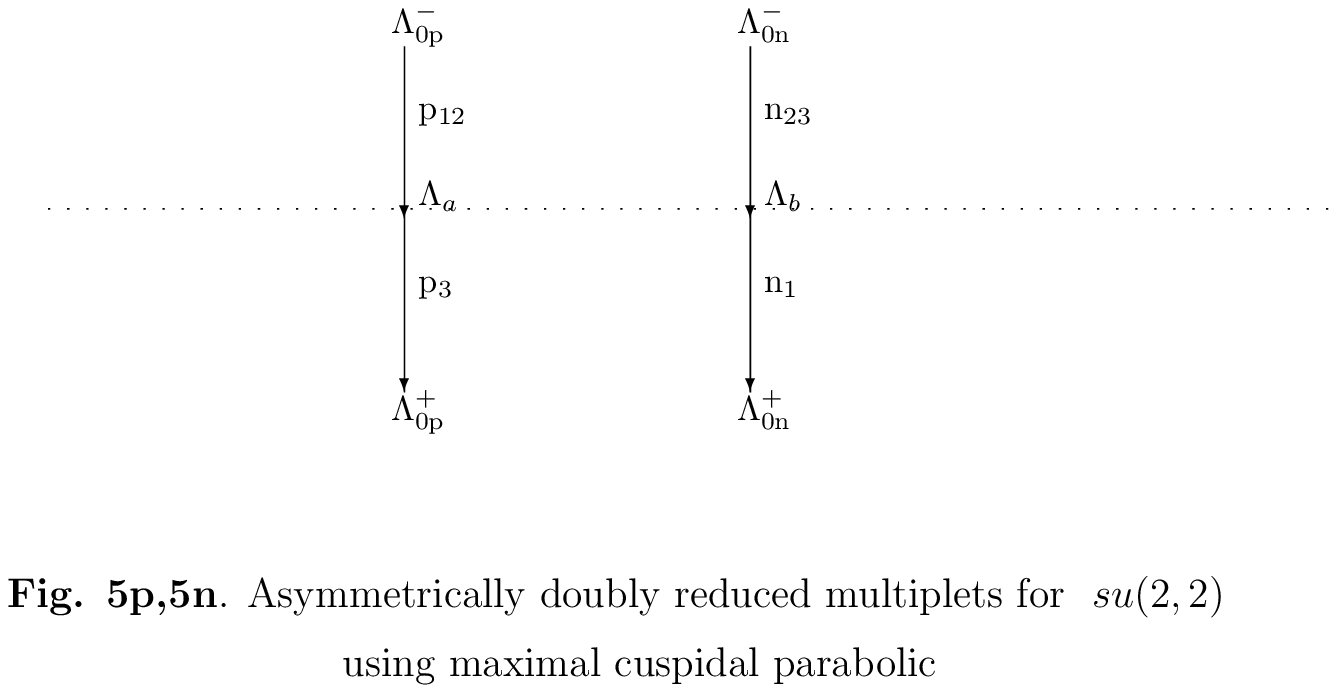}{16cm}


\begin{thebibliography}{99}




\bibitem{Dobinv}V.K. Dobrev, Invariant Differential Operators for Non-Compact
Lie Groups: Parabolic Subalgebras,
Rev. Math. Phys. {\bf 20} (2008) 407-449; hep-th/0702152.


\bibitem{Dobparab} V.K. Dobrev, Invariant Differential Operators for Non-Compact
Lie Algebras Parabolically Related to Conformal Lie Algebras,
J. High Energy Phys. 02 (2013) 015, arXiv:1208.0409.


 \bibitem{Lan}R.P. Langlands, {\it On the classification of irreducible
representations of real algebraic groups}, Math. Surveys and
Monographs, Vol. 31 (AMS, 1988), first as IAS Princeton preprint
(1973).

\bibitem{KnZu}A.W. Knapp and G.J. Zuckerman, ``Classification theorems
for representations of semisimple groups'',
 in: Lecture Notes in Math., Vol. 587 (Springer, Berlin,
1977) pp. 138-159; 
~Ann. Math. {\bf 116} (1982) 389-501.

 \bibitem{DMPPT}V.K. Dobrev, G. Mack, V.B. Petkova, S.G. Petrova and
I.T. Todorov,  {\it Harmonic Analysis on the $n$-Dimensional Lorentz
Group and Its Applications to Conformal Quantum Field Theory},
Lecture Notes in Physics, Vol. 63 (Springer, Berlin, 1977).

\bibitem{Knapp}A.W. Knapp, {\it Representation Theory of Semisimple
Groups (An Overview Based on Examples)}, (Princeton Univ. Press,
1986).


\bibitem{Dob}V.K. Dobrev,
`Canonical construction of intertwining differential operators
associated with representations of real semisimple Lie groups'',
Rept. Math. Phys. {\bf 25}, 159-181 (1988) ; first as ICTP Trieste
preprint IC/86/393 (1986).


 \bibitem{BGG}I.N. Bernstein, I.M. Gel'fand and S.I. Gel'fand,
``Structure of representations generated by highest weight vectors'',
Funkts. Anal. Prilozh. {\bf 5} (1) (1971) 1-9; English
translation: Funct. Anal. Appl. {\bf 5} (1971) 1-8.

\bibitem{Dix}J. Dixmier, {\it Enveloping Algebras}, (North Holland, New
York, 1977).

 \bibitem{Dobmul}V.K. Dobrev,
Multiplet classification of the reducible elementary
representations of real semi-{\allowbreak}simple Lie groups: the \
$SO_e(p,q)$ example,
 Lett. Math. Phys. {\bf 9}, 205-211 (1985).


\bibitem{DoJMP}V.K. Dobrev,
Elementary representations and intertwining operators for $SU(2,2)$: I, ~
J. Math. Phys. {\bf 26} (1985) 235-251.

\bibitem{HC}Harish-Chandra,
"Representations of semisimple Lie groups: IV,V",
Am. J. Math. {\bf 77} (1955) 743-777, {\bf 78} (1956) 1-41.

\bibitem{KnSt}A.W. Knapp and E.M. Stein,
``Intertwining operators for semisimple groups'',
Ann. Math. {\bf 93} (1971) 489-578; II : Inv. Math. {\bf 60} (1980) 9-84.



\bibitem{DoMo} V.K. Dobrev and P. Moylan,
Induced representations and invariant integral operators for $SU(2,2)$,
Fort. d. Physik, {\bf 42} (1994) 339-392.



\end{thebibliography}
\end{document}